\documentclass[11pt]{amsart}
\usepackage{amsmath}
\usepackage[pdftex, pdfstartview=FitH]{hyperref}

\theoremstyle{definition}

\newtheorem*{definition*}{Definition}

\theoremstyle{plain}
\newtheorem*{theorem*}{Theorem}
\newtheorem*{corollary*}{Corollary}
\newtheorem*{lemma*}{Lemma}
\newtheorem*{proposition*}{Proposition}

\newtheorem*{result*}{Main Result}

\numberwithin{equation}{section}

\newcommand{\Proof}{\noindent {\it Proof.\ }}

\newcommand{\dotdot}{{\, . \, . \,}}

\newcommand{\RP}{{\mathbb RP}}

\newcommand{\haken}{\mathbin{\hbox to 6pt{%
                 \vrule height0.4pt width5pt depth0pt
                 \kern-.37pt
                 \vrule height6pt width0.4pt depth0pt\hss}}}

\title{Asymmetry in Hilbert's Fourth Problem}

\author{J.C. \'Alvarez Paiva}
\address{J.C. \'Alvarez Paiva, U.M.R. CNRS 8524
         U.F.R. de Math\'ematiques, 59655 Villeneuve d'Ascq C\'edex, France}
\email{alvarez@math.univ-lille1.fr}

\dedicatory{To Fred Almgren, in memoriam}
\thanks{This work was partially supported by Programme ANR Blanc - SIMI 1 - ANR 2012}



\begin{document}%

\begin{abstract}
In the asymmetric setting, Hilbert's fourth problem asks to construct and study all (non-reversible) projective 
Finsler metrics: Finsler metrics defined on open, convex subsets of real projective $n$-space for which geodesics lie on projective lines. 
While asymmetric norms and Funk metrics provide many examples of essentially non-reversible projective metrics defined on {\it proper} convex
subsets of projective $n$-space, it is shown that {\it any projective Finsler metric defined on the whole projective $n$-space is the sum of a
reversible projective metric and an exact $1$-form.} The result admits a wide generalization to path geometries other than projective lines in
real projective spaces. 
\end{abstract}
\maketitle


\section{Introduction}
In the Summer of 1996 Eric Grinberg and Fred Almgren organized a very small and very informal meeting on Hilbert's fourth problem at the 
Institute of Advanced Studies. In the discussion after the talks Almgren asked whether I had considered asymmetric variational problems where 
the Lagrangian function does not satisfy $L(x,v)=L(x,-v)$ and, more precisely, whether every non-reversible Lagrangian on the plane 
whose extremals are straight lines is the sum of a reversible Lagrangian and the differential of a smooth function on the plane. I 
answered that any asymmetric norm in the plane provides an example that is not of this form and forgot about Almgren's question for the next 
fourteen years. In 2010 Gautier Berck surprised me with the following result:

\begin{theorem*}[Berck]
A projective Finsler metric defined on the whole real projective plane is the sum of a reversible projective Finsler metric and an exact $1$-form. 
\end{theorem*}

Berck's one-page proof is simple and elegant. However, it is also purely two-dimensional 
and applies only to projective Finsler metrics. The aim of this short note is to show through a different circle of ideas that 
the result also holds in higher dimensions and for many systems of paths other than projective lines. 

In order to state the main result of this note in a clear and succint form, let us recall that a Finsler metric is {\it geodesically
reversible} if every oriented geodesic can be reparametrized as a geodesic with the reverse orientation, and let us say that a Finsler 
metric is a {\it Zoll metric} if all of its geodesics are closed and of the same length.

\begin{result*}
Let $M$ be a manifold diffeomorphic to a compact rank-one symmetric space. If $L$ is a geodesically reversible Zoll metric on $M$, then
$L$ is the sum of a reversible Zoll metric and an exact $1$-form.  
\end{result*}

In the particular case of projective Finsler metrics on $\RP^n$ or $S^n$, this theorem completely solves the asymmetric version of
Hilbert's fourth problem by stating that, in fact, {\it there are no interesting non-reversible projective metrics on these spaces}. 
Notice however, that there are many interesting non-reversible projective metrics such as asymmetric norms or Funk metrics if we 
look at proper convex subsets of projective space.

\section{Preliminaries}

In this note a function $F: TM \rightarrow [0 \dotdot \infty)$ will be called a {\it Finsler metric} on the manifold $M$ if it
is smooth outside the zero section and its restriction to each tangent space $T_x M$ is a (possibly) asymmetric, quadratically 
convex norm $F_x$. Equivalently, the restriction of $F$ to each tangent space $T_x M$ is the support function of a quadratically
convex body $D_x^* M \subset T_x^* M$ that contains the origin in its interior. Finsler metrics for which all norms $F_x$, $x \in M$, 
are symmetric will be called {\it reversible Finsler metrics.} The reader may wonder why it is necessary to use two different terms---symmetry and 
reversibility---for the same concept, but this usage solves the problem of how to call the Finsler analogues of Riemannian 
symmetric spaces. Moreover, Finsler geometry is part of the geometry of Hamiltonian systems and the term reversible or
{\it time-reversible} Hamiltonian is well-established. 

The proof of our main result turns out to be surprisingly simple. Its key technical ingredients are the Brunn-Minkowski inequality 
and the Weinstein-Yang-Reznikov solution of what was then called the {\it weak Blaschke conjecture}. The key conceptual ingredient 
is the {\it Holmes-Thompson volume}.

\begin{definition*}
Let $(M,F)$ be an $n$-dimensional Finsler manifold and let $H : T^*M \rightarrow [0 \dotdot \infty)$ be the Hamiltonian function obtained from $F$ by 
the Legendre transform. The Holmes-Thompson volume of $(M,F)$ is the symplectic volume of the unit co-disc bundle $H \leq 1$ divided by the volume of 
the Euclidean unit ball of dimension $n$. 
\end{definition*} 

We may also define the Holmes-Thompson volume  of $(M,F)$ as the contact volume of the unit co-sphere bundle
$H = 1$ divided by the volume of the Euclidean unit ball of dimension $n$. Note that the Holmes-Thompson volume of a 
Riemannian manifold is its usual volume, which can also be defined as the Hausdorff measure of its associated metric. However,
for more general normed spaces and reversible Finsler manifolds the two notions of volume do not coincide. In (non-reversible) 
Finsler manifolds there does not seem to be a reasonable generalization of the Hausdorff measure. The reader is referred to~\cite{Thompson:1996} 
and~\cite{Alvarez-Thompson:2004} for more information about the Holmes-Thompson volume and other notions of volume on normed and Finsler spaces. 
For the purposes of this paper, we need only remark that the symplectic character of of the Holmes-Thompson volume makes it the volume of choice 
to trivially extend Riemannian results that only depend on the symplectic properties of the geodesic flow. One such result is the proof of the weak 
Blaschke conjecture by A.~Weinstein, C.T.~Yang, and A.~Reznikov (see \cite{Weinstein:1974}, \cite{Yang:1980, Yang:1991}, and \cite{Reznikov:1985,Reznikov:1994}). 
The following statement summarizes their results and puts them in a Finsler-geometric setting.

\begin{theorem*}[Weak Blaschke Conjecture]
Let $F$ and $G$ be two Zoll metrics on a manifold $M$ diffeomorphic to a compact rank-one symmetric space. If the lengths of 
the (prime) closed geodesics of $(M,F)$ and $(M,G)$ coincide, then their Holmes-Thompson volumes coincide as well.
\end{theorem*}

\section{Proof of the Main Theorem}

Let $a : TM \rightarrow TM$ be the map that sends each tangent vector $v$ to its opposite $-v$. The main result follows from a comparison between 
the properties of the Finsler metric $L$ and those of its symmetrization $\bar{L} := (L + L \circ a)/2$. First we show that symmetrization increases volume:

\begin{lemma*}
The Holmes-Thompson volume of a compact Finsler manifold $(M,L)$ is less than or equal to that of its symmetrization $(M,\bar{L})$.
Moreover, the volumes are equal if and only if $L$ equals its symmetrization modulo the addition of a $1$-form. In other words, 
$L = \bar{L} + \beta$, where $\beta$ is a $1$-form on $M$.   
\end{lemma*}

\Proof 
For each $x \in M$ denote by $D_x^* M \subset T^*M$ the convex body supported by the restriction of $L$ to $T_xM$. Note that the   
restriction of $\bar{L}$ to $T_xM$ is the support function of the symmetrized body $(D_x^* M + -D_x^* M)/2$. A well-known consequence
of the Brunn-Minkowski inequality (see, for example, \cite{Bonnesen-Fenchel}) is that the volume of a convex body $K$ is less than or equal
to that of its symmetrization $(K + -K)/2$ with equality if and only if $K$ is symmetric with respect to a (unique) point $p$ in its 
interior. Applying the theorem at every cotangent space $T_xM$, $x \in M$, proves the lemma. Note that the $1$-form will be smooth if
the bodies $D_x^*M$ depend smoothly on the base point, which is the case.
\qed

\begin{definition*}
Two Finsler metrics are {\it projectively equivalent} if they have the same geodesics up to orientation-preserving reparametrization. 
\end{definition*}

\begin{lemma*}
A geodesically reversible Finsler metric $L$ is projectively equivalent to its symmetrization $\bar{L} := (L + L\circ a)/2$.  
\end{lemma*}

\Proof
Remark only that projectively equivalent Finsler metrics form a convex cone and that $L$ is geodesically reversible if and only if $L$ and 
$L \circ a$ are projectively equivalent. \qed

\begin{lemma*}[Crampin~\cite{Crampin:2005}]
If $(M,F)$ is a reversible Finsler manifold and $\beta$ is a $1$-form on $M$, the Lagrangian $F + \beta$ is not geodesically reversible 
unless $\beta$ is closed. 
\end{lemma*}

\noindent {\it Proof of the Main Result.}
Since $L$ is geodesically reversible, it is projectively equivalent to its symmetrization $\bar{L}$ and, therefore, $L$ and $\bar{L}$ are two Zoll metrics 
such that the lengths of their (prime) closed geodesics coincide. By the solution of the Weak Blaschke Conjecture, the Holmes-Thompson volume of 
$(M,L)$ equals that of $(M,\bar{L})$. Applying the equality case of the Brunn-Minkowski inequality, we conclude that $L$ equals its symmetrization modulo the 
addition of a $1$-form $\beta$. Crampin's Lemma now implies that $\beta$ is closed and---given the topology of $M$---exact. \qed

\bibliography{../paperbib}
\bibliographystyle{amsplain}
\end{document}